\documentclass[a4paper,twoside]{article}

\newcommand{\heuteIst}{June 30, 2005 }

\usepackage{amsmath,amsthm,amsfonts,latexsym,amscd,amssymb,enumerate}
\usepackage[hypertex,backref]{hyperref}

\swapnumbers
\theoremstyle{plain}
\newtheorem{theorem}{Theorem}[section]
\newtheorem{lemma}[theorem]{Lemma}
\newtheorem{corollary}[theorem]{Corollary}
\newtheorem{proposition}[theorem]{Proposition}

\theoremstyle{definition}
\newtheorem{definition}[theorem]{Definition}
\newtheorem{example}[theorem]{Example}

\theoremstyle{remark}
\newtheorem{remark}[theorem]{Remark}

\newtheorem{exercise}[theorem]{Exercise}

\newcommand{\reals}{\mathbb{R}}
\newcommand{\complexs}{\mathbb{C}}
\newcommand{\naturals}{\mathbb{N}}
\newcommand{\integers}{\mathbb{Z}}

\DeclareMathOperator{\id}{id}

\newcommand{\abs}[1]{\left\lvert#1\right\rvert} 
\newcommand{\norm}[1]{\left\lVert#1\right\rVert}

\newcommand{\tensor}{\otimes}

\newcommand{\iso}{\cong}

\newcommand{\semiProd}{\rtimes}
\newcommand{\semiprod}{\semiProd}

\DeclareMathOperator{\im}{im}      
\DeclareMathOperator{\Hom}{Hom}    

\DeclareMathOperator{\diag}{diag}

\DeclareMathOperator{\tr}{tr}
\DeclareMathOperator{\pr}{pr}

\DeclareMathOperator{\Diffeo}{Diffeo}

\newcommand{\forget}[1]{}

\newcommand{\innerprod}[1]{\langle #1 \rangle}

{\catcode`@=11\global\let\c@equation=\c@theorem}

\allowdisplaybreaks[2]

\begin{document}
\pagestyle{myheadings}
\markboth{Thomas Schick}{Loop groups and string topology}

\date{Last compiled \today; last edited  \heuteIst or later}

\title{Loop groups and string topology\\ Lectures for the summer
  school algebraic groups\\
 G{\"o}ttingen, July 2005}

\author{Thomas Schick\thanks{
\protect\href{mailto:schick@uni-math.gwdg.de}{e-mail:
  schick@uni-math.gwdg.de}
\protect\\
\protect\href{http://www.uni-math.gwdg.de/schick}{www:~http://www.uni-math.gwdg.de/schick}
\protect\\
Fax: ++49 -551/39 2985
}\\
Uni G{\"o}ttingen\\
Germany}
\maketitle

\section{Introduction}
\label{sec:introduction}

Let $G$ be a compact Lie group. Then the space $LG$ of all maps from
the circle $S^1$ to $G$ becomes a group by pointwise
multiplication. Actually, there are different variants of $LG$,
depending on the classes of maps one considers, and the topology to be
put on the mapping space. In these lectures, we will always look at
the space of smooth (i.e.~$C^\infty$) maps, with the topology of
uniform convergence of all derivatives.

These groups certainly are not algebraic groups in the usual sense of
the word. Nevertheless, they share many properties of algebraic groups
(concerning e.g.~their representation theory). There are actually
analogous objects which are very algebraic (compare e.g.~\cite{MR1961134}), and
it turns out that those have properties remarkably close to those of the
smooth loop groups.

The lectures are organized as follows.
\begin{enumerate}
\item Lecture 1: Review of compact Lie groups and their representations, basics
of loop groups of compact Lie groups.
\item Lecture 2: Finer properties of loop groups
\item Lecture 3: the representations of loop groups (of positive energy)
\end{enumerate}

The lectures and these notes are mainly based on the excellent
monograph ``Loop groups'' by Pressley and Segal \cite{MR900587}.

\section{Basics about compact Lie groups}
\label{sec:basics-about-compact}

\begin{definition}
  A Lie group $G$ is a smooth manifold $G$ with a group structure,
  such that the map $G\times G\to G; (g,h)\mapsto gh^{-1}$ is smooth.

  The group acts on itself by left multiplication: $l_g(h)=gh$. A
  vector field $X\in\Gamma(TG)$ is called \emph{left invariant}, if
  $(l_g)_*X=X$ for each $g\in G$. The space of all left invariant
  vector fields is called the \emph{Lie algebra} $Lie(G)$. If we
  consider vector fields as derivations, then the commutator of two
  left invariant
  vector fields again is a left invariant vector field. This defines
  the Lie bracket $[\cdot,\cdot]\colon Lie(G)\times Lie(G)\to Lie(G); [X,Y]
  = XY-YX$.

  By left invariance, each left invariant vector field is determined
  uniquely by its value at $1\in G$, therefore we get the
  identification $T_1G\iso Lie(G)$; we will frequently use both
  variants.

  To each left invariant vector field $X$ we associate its flow
  $\Psi_X\colon G\times\reals\to G$ (a priori, it might only be
  defined on an open subset of $G\times \{0\}$). We define the
  \emph{exponential map} 
  \begin{equation*}
    \exp\colon Lie(G)\to G; X\mapsto \Psi_X(1,1).
  \end{equation*}
  This is defined on an open subset of $0\in G$. The differential
  $d_0\exp\colon Lie(G)\to T_1(G)=Lie(G)$ is the identity, therefore on a
  suitably small open neighborhood of $0$, $\exp$ is a diffeomorphism
  onto its image.

  A maximal torus $T$ of the compact Lie group $G$ is a Lie subgroup
  $T\subset G$ which is isomorphic to a torus $T^n$ (i.e.~a product of
  circles) and which has
  maximal rank among all such. It's a theorem that for a connected
  compact Lie group $G$ and a given maximal torus $T\subset G$, an
  arbitrary connected abelian Lie subgroup $A\subset G$ is conjugate to a
  subgroup of $T$.
\end{definition}

\begin{example}
  The group $U(n):=\{A\in M(n,\complexs)\mid AA^*=1\}$ is a Lie group,
  a Lie submanifold of the group $Gl(n,\complexs)$ of all invertible
  matrices.

  In this case, $T_1U(n)=\{A\in M(n,\complexs)\mid A+A^*=0\}$. The
  commutator of $T_1U(n)=L(U(n))$ is the usual commutator of matrices:
  $[A,B]=AB-BA$ for $A,B\in T_1U(n)$.

  The exponential map for the Lie group $U(n)$ is the usual exponential map of matrices, given
  by the power series:
  \begin{equation*}
  \exp\colon T_1 U(n)\to U(n); A\mapsto \exp(A)=\sum_{k=0}^\infty
  A^k/k!.
\end{equation*}
  The functional equation shows that the image indeed belongs to $U(n)$.

  Similarly, $Lie(SU(n))=\{A\in M(n,\complexs)\mid A+A^*=0, \tr(A)=0\}$.

\end{example}

\begin{theorem}
  If $G$ is a connected compact Lie group, then $\exp\colon Lie(G)\to G$ is surjective.
\end{theorem}

\begin{definition}\label{def:complexification}
  Every compact connected Lie group $G$ can be realized as a Lie-subgroup
  of $SU(n)$ for big enough $n$. It follows that its Lie algebra
  $Lie(G)$ is a sub-Lie algebra of $Lie(U(n))=\{A\in
  M(n,\complexs)\mid A^*=-A\}$. Therefore, the complexification
  $Lie(G)\tensor_\reals\complexs$ is a sub Lie algebra of
  $Lie(U(n))\tensor_\reals\complexs = \{A+iB\mid A,B\in
  M(n,\complexs), A^*=-A, (iB)^*=-(iB)\}= M(n,\complexs)$. Note that
  the bracket (given by the commutator) is complex linear on these
  complex vector spaces.

  The corresponding sub Lie group $G_\complexs$ of $Gl(n,\complexs)$ (the simply
  connected Lie group with Lie algebra $M(n,\complexs)$) with Lie
  algebra $Lie(G)\tensor_\reals\complexs$ is called the
  complexification of $G$.

  $G_\complexs$ is a \emph{complex Lie group}, i.e.~the manifold $G_\complexs$ has a
  natural structure of a complex manifold (charts with holomorphic
  transition maps), and the composition $G_\complexs\times
  G_\complexs\to G_\complexs$ is holomorphic.
\end{definition}

\begin{exercise}
  Check the assertions made in Definition \ref{def:complexification},
  in particular about the complex structure.
\end{exercise}

\begin{example}
  The complexified Lie algebra $Lie(SU(n))\tensor_\reals\complexs=\{
  A\in M(n,\complexs)\mid \tr(A)=0\}$, and the complexification of
  $SU(n)$ is $Sl(n,\complexs)$.
\end{example}

\begin{definition}
  Let $G$ be a compact Lie group. It acts on itself by conjugation:
  $G\times G\to G; (g,h)\mapsto ghg^{-1}$.

  For fixed $g\in G$, we can take the differential of the
  corresponding map $h\mapsto ghg^{-1}$ at $h=1$. This defines the
  adjoint representation $ad\colon G\to Gl(Lie(G))$.

  We now decompose $Lie(G)$ into irreducible sub-representations for
  this action, $Lie(G)=\frak{g}_1\oplus \cdots\oplus \frak{g}_k$. Each
  of these are Lie subalgebras, and we have
  $[\frak{g}_i,\frak{g}_j]=0$ if $i\ne j$.

  $G$ is called {semi-simple} if non of the summands is one-dimensional,
  and \emph{simple} if there is only one summand, which additionally
  is required  not to be one-dimensional.
\end{definition}

\begin{remark}
  The simply connected simple compact Lie groups have been classified,
  they consist of $SU(n)$, $SO(n)$, the symplectic groups $Sp_n$ and
  five exceptional groups (called $G_2$, $F_4$, $E_6$, $E_7$, $E_8$). 
\end{remark}

\begin{definition}
  Let $G$ be a compact Lie group with a maximal torus $T$. $G$ acts
  (induced from conjugation) on $Lie(G)$ via the adjoint
  representation, which induces a representation on
  $Lie(G)\tensor_\reals\complexs=: \mathfrak{g}_\complexs$. This Lie
  algebra contains the Lie algebra $\mathfrak{t}_\complexs$ of the
  maximal torus, on which $T$ acts trivially (since $T$ is abelian).

  Since $T$ is a \emph{maximal} torus, it acts non-trivial on every
  non-zero vector of the complement.

  As every finite dimensional representation of a torus, the complement decomposes into a
  direct sum $\bigoplus \mathfrak{g}_\alpha$, where on each
  $\mathfrak{g}_\alpha$, $t\in T$ acts via multiplication with
  $\alpha(t)\in S^1$, where $\alpha\colon T\to S^1$ is a homomorphism,
  called the \emph{weight} of the summand $\mathfrak{g}_\alpha$.

  We can translate the homomorphisms $\alpha\colon T\to S^1$ into
  their derivative at the identity, thus getting a linear map
  $\alpha'\colon \mathfrak{t}\to\reals$, i.e.~an element of the dual space
  $\mathfrak{t}^*$, they are related by
  $\alpha(\exp(x))=e^{i\alpha'(x)}$.

  This way, we think of the group of characters $\hat T=\Hom(T,S^1)$
  as a lattice in $t^*$, called the lattice of weights. It contains
  the set of \emph{roots}, i.e.~the non-zero weights occurring in the
  adjoint representation of $G$.
\end{definition}

\begin{remark}\label{1_dim} 
  As $\mathfrak{g}_\complexs$ is the complexificatoin of a real representation,
if $\alpha$ is a root of $G$, so is $- \alpha$, with
$\mathfrak{g_{-\alpha}}=\overline{\mathfrak{g_\alpha}}$.
 
 It is a theorem that the subspaces $\mathfrak{g_\alpha}$ are always $1$-dimensional.
\end{remark}

\begin{example}
  For $U(n)$, $Lie(U(n))_\complexs=M(n,\complexs)$. A maximal torus is given by
  the diagonal matrices (with diagonal entries in $S^1$). 
  We can then index the roots by pairs $(i,j)$ with $1\le i\ne j\le n$. We have
$\alpha_{ij}(\diag(z_1,\dots,z_n))=z_iz_j^{-1}\in S^1$. The corresponding
subspace $\mathfrak{g}_{ij}$ consists of matrices which are zero except for
the $(i,j)$-entry. As an element of $\mathfrak{t}^*$ it is given by the linear
map $\reals^n\to\reals; (x_1,\dots,x_n)\mapsto x_i-x_j$.
\end{example}

\begin{definition}
  According to Remark \ref{1_dim}, each of the spaces $\mathfrak{g}_\alpha$ is
$1$-dimensional. Choose a vector $0\ne e_\alpha\in\mathfrak{g}_\alpha$, such
that $e_{-\alpha}=\overline{e_\alpha}$. Then
$h_\alpha:=-i[e_\alpha,e_{-\alpha}]\in\mathfrak{t}$ is non-zero. We can
normalize the vector $e_\alpha$ in such a way that
$[h_\alpha,e_\alpha]=2ie_\alpha$. Then $h_\alpha$ is canonically determined by
$\alpha$, it is called the \emph{coroot} associated to $\alpha$.

We call $G$ a \emph{simply knitted} Lie group, if there is a
$G$-invariant inner 
product on $Lie(G)$ such that $\innerprod{h_\alpha,h_\alpha}=2$
 for all roots $\alpha$.
 
 This inner product gives rise to an isomorphism
$\mathfrak{t}\iso\mathfrak{t}^*$; this isomorphism maps $h_\alpha$ to $\alpha$.
 \end{definition}
 
 \begin{example}
   For $SU(2)$, $T=\{\diag(z,z^{-1})\mid z\in S^1\}$, there is only one pair
of roots $\alpha,-\alpha$, obtained by restriction of the roots $\alpha_{12}$
 and $\alpha_{21}$ of $U(2)$ to the (smaller) maximal torus of $SU(2)$. We
get $\alpha(\diag(z,z^{-1}))=z^2$,  $e_\alpha=\left(\begin{smallmatrix} 0&1\\ 0 &0
\end{smallmatrix}\right)$, $h_\alpha=\left(\begin{smallmatrix} i & 0\\ 0 &
-i\end{smallmatrix}\right)$.

All the groups $U(n)$ and $SU(n)$ are simply knitted, as well as $SO(2n)$. The
canonical inner product of $Lie(U(n))_\complexs=M(n,\complexs)$ is given by
$\innerprod{A,B}=-tr(A^*B)/2$.
 \end{example}

 \begin{proposition}\label{prop:generation by Lie su2}
   Given any semi-simple Lie group, the Lie algebra is spanned by the
   vectors $e_\alpha,e_{-\alpha},h_\alpha$ as $\alpha$ varies over the
   set of (positive) roots.

 For each such triple one can define a Lie algebra homomorphism
 $Lie(SU(2))\to Lie(G)$ which maps the standard generators of
 $Lie(SU(2))$ to the chosen generators of $Lie(G)$.

 We get the technical useful result that for a semi-simple compact Lie
 group $G$ the Lie algebra $Lie(G)$ is generated by the images of
 finitely many Lie algebra homomorphisms from $Lie(SU(2))$ to $L(G)$.
 \end{proposition}

\begin{definition}
  Let $G$ be a compact Lie group with maximal torus $T$. Let $N(T):=\{g\in
G\mid gTg^{-1}\subset T\}$ be the normalizer of $T$ in $G$. The \emph{Weyl
group} $W(T):= N(T)/T$ acts on $T$ by conjugation, and via the adjoint
representation also on $\mathfrak{t}$ and then on $\mathfrak{t}^*$. It is a
finite group.

This action preserves the lattice of weights $\hat T\subset \mathfrak{t}^*$,
and also the set of roots.

Given any root $\alpha$ of $G$, the Weyl group contains an element $s_\alpha$
of order $2$. It acts on $\mathfrak{t}$ by reflection on the hyperplane
$H_\alpha:=\{X\mid\alpha(X)=0\}$. More precisely we have
$s_\alpha(X)=X-\alpha(X)h_\alpha$, where $h_\alpha\in\mathfrak{t}$ is the
coroot associated to the root $\alpha$.

The reflections $s_\alpha$ together generate the Weyl group $W(G)$.

The Lie algebra $\mathfrak{t}$ is decomposed into the union of the root
hyperplanes $H_\alpha$ and into their complement, called the set of
\emph{regular} elements. The complement decomposes into finitely many
connected components, which are called the \emph{Weyl chambers}. One chooses
one of these and calls it the \emph{positive Weyl chamber}. The Weyl group
acts freely and transitively on the Weyl chambers.

A root of $G$ is called positive or negative, if it assumes positive or
negative values on the positive chamber. A positive root $\alpha$ is called
\emph{simple}, if its hyperplane $H_\alpha$ is a wall of the positive chamber.
\end{definition}

\begin{example}
  For the group $U(n)$, the Weyl group is the symmetric group $S_n$, it acts
by permutation of the diagonal entries on the maximal torus. Lifts of the
elements of $W(U(n))$ to $N(T)\subset U(n)$ are given by the permutation
matrices.

For $SU(n)$, the Weyl group is the alternating group $S_n$, again acting by
permutation of the diagonal entries of the maximal torus.

Since these groups are simply knitted, we can depict their roots,
coroots etc.~in a simply Euclidean picture of
$\mathfrak{t}\iso\mathfrak{t}^*$. We look at the case $G=SU(3)$ where
$Lie(T)$ is $2$-dimensional. 

In the lecture, a picture is drawn.

In general, the positive roots of $U(n)$ or $SU(n)$ can be chosen to
be the roots $\alpha_{ij}$ with $i<j$. 
\end{example}

\begin{theorem}
  Let $G$ be a compact semi-simple Lie group. There is a one-to-one
  correspondence between the irreducible representations of $G$ and
  the set of dominant weights, where a weight $\alpha$ is called
  dominant if $\alpha(h_\beta)\ge 0$ for each positive root $\beta$ of
  $G$.
\end{theorem}
\begin{proof}
  We don't proof this theorem here; we just point out that the
  representation associated to a dominant weight is given in as the
  space of holomorphic sections of a line bundle over a homogeneous
  space of $G_\complexs$ associated to the weight.
\end{proof}

\begin{remark}
  One of the goals of these lectures is to explain how this results
  extends to loop groups.
\end{remark}

\section{Basics about loop groups}
\label{sec:basics-about-loop}

\begin{definition}
  An infinite dimensional smooth manifold (modeled on a locally
  convex complete
  topologically vector space $X$) is a topological space $M$ together
  with a collection of charts $x_i\colon U_i\to V_i$, with open subset
  $U_i$ of $M$ and $V_i$ of $X$ and a homeomorphism $x_i$, such that
  the change of coordinate maps $x_j\circ x_i^{-1}\colon V_i\to V_j$
  are smooth maps between the topological vector space $X$.

  Differentiability of a map $f\colon X\to X$ is defined in terms of
  convergence of difference quotients, if it exists, the differential
  is then a map 
  \begin{equation*}Df\colon X\times X\to X; (v,w)\mapsto \lim_{t\to 0} \frac{f(u+tv)-f(u)}{t}\end{equation*}
   (where the second variable encodes
  the direction of differentiation), iterating this, we define higher
  derivatives and the concept of smooth, i.e.~$C^\infty$ maps.
\end{definition}

\begin{definition} Let $G$ be a compact Lie group with Lie algebra
  $Lie(G)=T_1G$ e( a finite dimensional vector space).
  For us, the model topological vector space $X$ will be
  $X=C^\infty(S^1,Lie(G))$. Its topology is defined by the collection of
  semi-norms $q_i(f):=\sup_{x\in S^1}\{\abs{\partial^{i}\phi/\partial
    t^k (x)}\}$ for any norm on $Lie(G)$.

  This way, $X$ is a complete separable (i.e.~with a countable dense
  subset)  metrizable topological vector space. A sequence of smooth
  functions $\phi_k\colon S^1\to Lie(G)$ converges to $\phi\colon S^1\to
  Lie(G)$ if and only if the functions and all their derivatives
  converge uniformly.
\end{definition}

\begin{lemma}\label{lem:LG_is_smooth}
There is a canonical structure of a smooth infinite dimensional
manifold on $LG$. A chart around the constant loop $1$ is given by the
exponential map $C^\infty(S^1,U)\to LG; \chi\mapsto \exp\circ
\chi$ where $U$ is a sufficiently small open neighborhood of $0\in
T_1(G)$. Charts around any other point $f\in LG$ are obtained by
translation with $f$, using the group structure on $LG$.

We define the topology on $LG$ to be the finest topology (as many open
sets as possible) such that all the above maps are continuous.

It is not hard to see that the transition functions are then actually
smooth, and that the group operations (defined pointwise) are smooth maps.
\end{lemma}

\begin{exercise}
  Work out the details of Lemma \ref{lem:LG_is_smooth}, and show that
  it extends to the case where $S^1$ is replaced by any compact smooth
  manifold $N$.
\end{exercise}

\begin{remark}
  There are many variants of the manifold $LG$ of loops in $G$. Quite
  useful are versions which are based on maps of a given Sobolev
  degree, one of the advantages being that the manifolds are then
  locally Hilbert spaces. 
\end{remark}

\begin{exercise}
  Show that for $G=SU(2)$ the group $LG$ is connected, but the
  exponential map is not surjective 
\end{exercise}

\begin{exercise}
  There are other interesting infinite dimensional Lie groups. One
  which is of some interest for loop groups is $Diffeo(S^1)$, the group
  of all diffeomorphisms of $S^1$ (and also its identity component of
  orientation preserving diffeomorphism).

  Show that this is indeed a Lie group, with Lie algebra (and local
  model for the smooth structure) $Vect(S^1)$, the space of all smooth
  vector fields. $\exp\colon Vect(S^1)\to Diffeo(S^1)$ maps a vector
  field to the (time 1) flow generated by it.

  Show that there are no neighborhoods $U\subset Vect(S^1)$ of $0$ and
  $V\subset Diffeo(S^1)$ of $\id_{S^1}$ such that $\exp|\colon U\to V$
  is injective or surjective.
\end{exercise}

\begin{lemma}\label{lem:decompose}
 Consider the subgroup of based loops $\Omega(G)=\{f\colon S^1\to G\in
 LG\mid f(1)=1\}\subset LG$, and the subgroup of constant loops
 $G\subset LG$.

 The multiplication map $G\times \Omega(G)\to LG$ is a diffeomorphism.
\end{lemma}
\begin{proof}
  The inverse map is given by $LG\to G\times \Omega(G); f\mapsto (f(1),
  f(1)^{-1}f)$. Clearly the two maps are continuous and inverse to
  each other. The differentiable structure of $\Omega(G)$ makes the
  maps smooth.
\end{proof}

\begin{corollary}\label{corol:homotopy groups}
  It follows that the homotopy groups of $LG$ are easy to compute (in
  terms of those of $G$); the maps of Lemma \ref{lem:decompose} give
  an isomorphism $\pi_k(LG)\iso \pi_k(G)\oplus \pi_{k-1}(G)$.

  In particular, $LG$ is connected if and only if $G$ is connected and
  simply connected, else the connected components of $G$ are
  parameterized by $\pi_0(G)\times \pi_1(G)$.
\end{corollary}

\begin{definition}
  Embed the compact Lie group $G$ into $U(n)$.
  Using Fourier decomposition, we can now write the elements of $LG$ in
  the form $f(z)=\sum A_kz^k$ with $A_k\in M(n,\complexs)$.

  We define now a number of subgroups of $LG$.
  \begin{enumerate}
  \item $L_{pol}G$ consists of those loops with only finitely many
    non-zero Fourier coefficients. It is the union (over
    $N\in\naturals$) of the subsets of
    functions of degree $\le N$. The latter ones are compact subsets
    of $M(n,\complexs)^N$, and we give $L_{pol}G$ the direct limit topology.
  \item $L_{rat}G$ consists of those loops which are rational
    functions $f(z)$ (without poles on $\{\abs{z}=1\}$).
  \item $L_{an}G$ are those loops where the series $\sum A_kz^k$
    converges for some annulus $r\le \abs{z}\le 1/r$ with $0<r<1$. For
    fixed $0<r<1$ this is Banach Lie group (of holomorphic functions
    on the corresponding annulus), with the topology of uniform
    convergence. We give $L_{an}G$ the direct limit topology.
  \end{enumerate}
\end{definition}

\begin{exercise}
  In general, $L_{pol}G$ is not dense in $G$. Show this for the case $G=S^1$.
\end{exercise}

\begin{proposition}
  If $G$ is semi-simple, then $L_{pol}G$ is dense in $LG$.
\end{proposition}
\begin{proof}
   Set $H:=\overline{L_{pol}}G$, and $V:=\{X\in
   C^\infty(S^1,Lie(G)\mid \exp(tX)\in H\forall t\in\reals\})\subset
   C^\infty(S^1,Lie(G))$. Then $V$ is a vector space, as 
   \begin{equation*}
\exp(X+Y) = \lim_{k\to\infty} \left(\exp(X/k)\exp(Y/k)\right)^k
\end{equation*}
converges in the $C^\infty$-topology if $X$ and $Y$ are close enough
to $0$. Since $\exp$ is continuous, $V$ is a closed subspace. It
remains to check that $V$ is dense in $L\mathfrak{g}$.

Because of Proposition \ref{prop:generation by Lie su2} it suffices to
check the statement for $SU(2)$. Here, we first note that the elements
$X_n\colon z\mapsto \left(
  \begin{smallmatrix}
    0 & z^n\\ -z^{-n} & 0
  \end{smallmatrix}\right)$ and $Y_n\colon z\mapsto \left(
  \begin{smallmatrix}
    0 & iz^n\\ iz^{-n} & 0
  \end{smallmatrix}\right)$ belong to $V$, since $X_n^2=Y_n^2=-1$, so
that $\exp(tX_n)=\sum_k (tX_n)^k/k!$ actually is a family of
polynomial loops. By linearity and the fact that $V$ is closed, every
loop of the form $z\mapsto \left(
  \begin{smallmatrix}
    0 & f(z)+ig(z)\\ -f(z)+ig(z) & 0
\end{smallmatrix}
\right)$ belongs to $V$, for arbitrary smooth real valued functions
$f,g\colon S^1\to \reals$ (use their Fourier decomposition). Since $V$
is finally invariant under conjugation by polynomial loops in $LSU(2)$, it
is all of $\mathfrak{su}(2)$.
\end{proof}

\subsection{Abelian subgroups of $LG$}
\label{sec:abelian-subgroups-lg}

We have already used the fact that each compact Lie group $G$ as a maximal
torus $T$, which is unique upto conjugation.

\begin{proposition}\label{prop:abelian subgroups}
 Let $G$ be a connected compact Lie group with maximal torus $T$.
  Let $A\subset LG$ be a (maximal) abelian subgroup. Then
  $A(z):=\{f(z)\in G\mit f\in A\}\subset G$ is a (maximal) abelian
  subgroup of $G$.

  In particular, we get a the maximal abelian subgroup $A_\lambda$ for each
  smooth map $\lambda\colon S^1\to \{T\subset 
  G\mid T \text{ maximal torus}\}\iso G/N$, where $N$ is the
  normalizer of $T$ in $G$, with $A_\lambda:=\{f\in LG\mid f(z)\in \lambda(z)\}$.

  The conjugacy class of $A_\lambda$ depends only on the homotopy
  class of $\lambda$. Since $\pi_1(G/N)=W=N/T$, $[S^1,G/N]$ is the set
  of conjugacy classes of elements of the Weyl group $W$.

  An element $w\in W$ acts on $T$ by conjugation, and the
  corresponding $A_\lambda$ is isomorphic to $\{\gamma\colon \reals\to
  T\mid \gamma(t+2\pi)=w^{-1}\gamma(t)w$ for all $t\in\reals$.
\end{proposition}
\begin{proof}
  It is clear that all evaluation maps of $A$ have abelian image. If
  all the image sets are maximal abelian, then $A$ is  maximal
  abelian.

  Since all maximal tori are conjugate, the action of $G$ on the set
  of all maximal tori is transitive, with stabilizer (by definition)
  the normalizer $N$, so that this space is isomorphic to $G/N$.

  Next, $W=N/T$ acts freely on $G/T$, with quotient $G/N$. Since $G/T$
  is simply connected (a fact true for every connected Lie group), we
  conclude from covering theory that
  $\pi_1(G/T)\iso W$. But then the homotopy classes of non base-point
  preserving maps are bijective to the conjugacy classes of the
  fundamental group.

  Given a (homotopy class of maps) $\lambda\colon S^1\to G/T$,
  represented by $w\in W=N/T$ (with a lift $w'\in N$ and with $x\in
  Lie(G)$ such that $\exp(2\pi x)=w'$), we can choose
  $\lambda$ with $\lambda(z)= \exp(zx) T \exp(-zx)$ (recall that the
  bijection between $G/N$ and the set of maximal tori is given by
  conjugation of $T$). 

  We then get a bijection from $A_\lambda$ to the twisted loop group
  of the assertion by sending $f\in A_\lambda$ to $\tilde
  f\colon\reals\to T: t\mapsto \exp(-tx)f(t)\exp(tx)$.
\end{proof}

\begin{example}
  The most obvious maximal abelian subgroup of a compact Lie group $G$
  with maximal torus $T$ is $LT$, which itself contains in particular
  $T$ (as subgroup of constant loops).
\end{example}

\begin{exercise}
  Find other maximal abelian subgroups of $LG$.
\end{exercise}

\begin{example}
  If $G=U(n)$ with maximal torus $T$, its Weyl group $W$ is isomorphic to the symmetric group
  $S_n$. Given a cycle $w\in W=S_n$, the corresponding maximal abelian
  subgroup $A_w\subset LU(n)$ is isomorphic to $LS^1$.

  More generally, if $w$ is a product of $k$ cycles (possibly of
  length $1$), then $A_w$ is a
  product of $k$ copies of $LS^1$.
\end{example}
\begin{proof}
   We prove the statement if $w$ consists of one cycle (of length $n$). Then, in
   the description of Proposition \ref{prop:abelian subgroups}, $A_w$
   consists of functions $\reals\to T$ which are periodic of period
   $2\pi l$. Moreover, the different components are all determined by
   the first one, and differ only by a translation in the argument of
   $2\pi$ or a multiple of $2\pi$.
\end{proof}

\begin{proposition}
  If $G$ is semi-simple and compact, then $LG_0$, the component of the
  identity of $LG$, is perfect.
\end{proposition}

\section{Root system and Weyl group of loop groups}
\label{sec:root-system-weyl}

\begin{definition}
  Let $G$ a compact Lie group with maximal torus $T$. Consider its complexified Lie algebra
  $Lie(LG)_\complexs= L\mathfrak{g}_\complexs$.

  It carries the action of $S^1$ by reparametrization of loops:
  $(z_0\cdot X)(z):= X(z_0z)$.

  We get a corresponding decomposition of $L\mathfrak{g}_\complexs=
  \bigoplus_{k\in \integers}\mathfrak{g}_{\complexs} z^k$   into
  Fourier components (the sums
  has of course to be completed appropriately).

  The action of $S^1$ used in the above decomposition still commutes
  with the adjoint representation of the subgroup $T$ of constant
  loops, so the summands can further be decomposed according to the
  weights of the action of $T$, to give a decomposition
  \begin{equation*}
     L\mathfrak{g}_\complexs = \oplus_{(k,\alpha)\in \integers\times
       \hat T} \mathfrak{g}_{(k,\alpha)} z^k.
   \end{equation*}
   The index set $\integers\times \hat T$ is the Pontryagin dual of
   $S^1\times T$ (i.e.~the set of all homomorphisms $S^1\times T \to
   S^1$). Those homomorphisms which occur (now also with possibly
   $\alpha=0$) are called the \emph{roots} of $LG$.

   The (infinite) set of roots of $LG$ is permuted by the so called
   \emph{affine Weyl group}
   \begin{equation*}
 W_{aff}=N(T\times S^1)/(T\times S^1),
\end{equation*}
   considered inside the semidirect product $LG\semiProd S^1$, where
   we use the action of $S^1$ by reparametrization (rotation) of loops
   on $LG$ to construct the semidirect product. This follows because
   we decompose $L\mathfrak{g}_\complexs$ as a representation of
   $LG\semiProd S^1$ with respect to the subgroup $T\times S^1$..
\end{definition}

\begin{proposition}
  $W_{aff}$ is a semidirect product $\Hom(S^1,T)\semiProd W$, where
  $W$ is the Weyl group of $G$, with its usual action on the target $T$.
\end{proposition}
\begin{proof}
  Clearly $\Hom(S^1,T)$ is a subgroup of $LG$ which conjugates the
  constant loops with values in $T$
  into itself, and the action of the constant loops with values in the
  normalizer $N$ does the same (and factors through $W$). It is also
  evident that all of $S^1$ belongs to the normalizer of $T\times S^1$.

  On the other hand, if $R_{z_0}\in S^1\in LG\semiProd S^1$ acts by
  rotation by $z_0\in S^1$, then for $f\in LG$ we get $f^{-1} R_{z-0}f=
  f^{-1}(\cdot)f(\cdot z_0)  R_{z_0} $.

  This belongs to $T\times S^1$ if and only if $z\mapsto
  f(z)^{-1}f(zz_0)$ is a constant element of $T$ for each $z_0$, i.e.~if and
  only if $z\mapsto f(1)^{-1}f(z)$ is a group homomorphism $z\to T$.
 Additionally, $f$ conjugates $T$ to itself if and only if $f(1)\in
 N$. Therefore, $N(T\times S^1)$ is the product of $N$,
 $\Hom(S^1,T)$ and $S^1$ inside $LG\semiProd S^1$, and the quotient by
 $T\times S^1$ is as claimed.
\end{proof}

\begin{definition}
  We think of the weights and roots of $LG$ not as linear forms on
  $\reals\times \mathfrak{t}$ (derivatives of elements in
  $\Hom(S^1\times T,S^1)$), but rather as \emph{affine linear}
  functions on $\mathfrak{t}$, where we identify $\mathfrak{t}$ with
  the hyperplane $1\times \mathfrak{t}$ in $\reals\times
  \mathfrak{t}$; this explains the notation ``affine roots''.

  Moreover, the group $W_{aff}$ acts linearly on
  $\reals\times\mathfrak{t}$ and preserves $1\times\mathfrak{t}$,
  where $\lambda\in \Hom(S^1,T)$ acts by translation by
  $\lambda'(1)\in \mathfrak{t}$.

  An affine root $(k,\alpha)$ is (for $\alpha\ne 0$) determined (upto
  sign) by the affine hyperplane $H_{k,\alpha}:=\{x\in\mathfrak{t}\mid
  \alpha(x)=-k\}\subset \mathfrak{t}$ which is the set where it
  vanishes.

  The collection of these hyperplanes is called the diagram of
  $LG$. It contains as the subset consisting of the $H_{0,\alpha}$ the
  diagram of $G$

  Recall that the connected components of the complement of all $H_{0,\alpha}$ were
    called the chambers of $G$, and we choose one which we call the
    positive chamber. The components of the complement of the diagram
    of $LG$ are called \emph{alcoves}. Each chamber contains a unique
    alcove which touches the origin, and this way the positive chamber
    defines a positive alcove, the set $\{x\in \mathcal{t}\mid
    0<\alpha(x)<1\text{ for all positive roots }\alpha\}$. An affine
    root is called positive or negative, if it has positive or
    negative values at the positive alcove. The positive affine roots
    corresponding to the walls of the positive alcove are called the
    simple affine roots.
\end{definition}

\begin{example}
  The diagram for $SU(3)$ is the tessellation of the plane by
  equilateral triangles.

  In general, if $G$ is a simple group, then each chamber is a
  simplicial cone, bounded by the $l$ hyperplanes
  $H_{0,\alpha_1},\dots,H_{0,\alpha_l}$ (where
  $\alpha_1,\dots,\alpha_l$ are the simple roots of $G$). There is a
  highest root $\alpha_{l+1}$ of $G$, which dominates all other roots
  (on the positive chamber). The positive alcove is then an
  $l$-dimensional simplex cut out of the positive chamber by the wall
  $H_{1,-\alpha_{l+1}}$, and we get $l+1$ simple affine roots of $LG$,
  $(0,\alpha_1),\dots, (0,\alpha_l),(1,-\alpha_{l+1})$.

  In general, if $G$ is semi-simple with $q$ simple factors, the positive alcove is a product
  of $q$ simplices, bounded by the walls of the simple roots of $G$,
  and walls $H_{1,-\alpha_i}$, $i=l+1,\dots,l+q$ being the highest
  weights of the irreducible summands of the adjoint representation. 
\end{example}

\begin{proposition}
  Let $G$ be a connected and simply connected compact Lie group. Then $W_{aff}$ is
  generated by reflections in the hyperplanes (the reflections in the
  walls of the positive alcove suffice), and it acts freely
  and transitively on the set of alcoves.
\end{proposition}
\begin{proof}
  From the theory of compact Lie groups, we know that $W$ is generated
  by the reflections at the $H_{0,\alpha}$, and that
  $\Hom(S^1,\mathfrak{t})$ is generated by the coroots $h_\alpha$.

  Recall that the reflection $s_\alpha$ was given by
  $s_\alpha(x)=x-\alpha(x)h_\alpha$. Recall that we normalized such
  that $\alpha(h_\alpha)=2$, therefore $-kh_\alpha/2\in H_{k,\alpha}$.

  Now the reflection $s_{k,\alpha}$ in the hyperplane $H_{k,\alpha}$
  is given by $s_{k,\alpha}(x)= x+kh_\alpha/2 -
  \alpha(x+kh_\alpha/2)h_\alpha - kh_\alpha/2 =
  s_\alpha(x)-kh_\alpha$. Since $s_\alpha\in W$ and $h_\alpha$ in
  $\Hom(S^1,T)$ (identified with the value of its derivative at $1$),
  $s_{k,\alpha}\in W_{aff}= \Hom(S^1,T)\semiProd W$.

  To show that these reflections generate $W_{aff}$, it suffices to
  show that they generate the translation by $h_\alpha$. But this is
  given by $s_{a,-\alpha}s_{0,\alpha}$.

  We now show that $W_{aff}$ acts transitively on the set of
  alcoves. For an arbitrary alcove $A$, we have to find $\gamma\in
  W_{aff}$ such that $\gamma A$ is the positive alcove $C_0$. Now the
  orbit $W_{aff}a$ of a point $p\in A$ is a subset $S$ of
  $\mathfrak{t}$ without accumulation points. Choose a point $c\in
  C_0$ and one of the points $b\in S$ with minimal distance to $c$. If
  $b$ would not belong to $C_0$, then $b$ and $c$ are separated by a
  wall of $C_0$, in which we can reflect $b$ to obtain another point
  of $S$, necessarily closer to $c$ than $b$.

  Since $W$ acts freely on the set of chambers, $W_{aff}$ acts freely
  on the set of alcoves (since each element of $W$ preserves the
  distance to the origin, and each translation moves the positive
  alcove away from the origin, only elements of $W$ could stabilize
  the positive alcove).
\end{proof}

\section{Central extensions of $LG$}
\label{sec:centr-extens-lg}

We want to study the representations of $LG$ for a compact Lie group
$G$. It turns out, however, that most of the relevant representations
are no honest representations but only projective representations.,
i.e.~$U_fU_g=c(f,g) U_{fg}$ for every $f,g\in LG$, where $U_f$ is the
operator by which $f$ acts, and $c(f,g)\in S^1$ is a scalar valued function
(a cocycle).

More precisely, the actions we consider are actions of central
extensions of $LG$ (in some sense defined by this cocycle). We don't
want to go into the details of the construction and classification of
these central extensions, but only state the main results.

Let $G$ be a compact connected Lie group.
\begin{enumerate}
\item $LG$ has many central extensions $1\to S^1\to \tilde LG\to LG
  \to 1$.
\item The corresponding Lie algebras $\tilde L\mathfrak{g}$ are
  classified as follows:
  for every symmetric invariant bilinear form
  $\innerprod{\cdot,\cdot}\colon \mathfrak{g}\times\mathfrak{g}\to
  \reals$ we get a form
  \begin{equation*}
    \omega\colon L\mathfrak{g}\times L\mathfrak{g}\to \reals;
    (X,Y)\mapsto \frac{1}{2\pi} \int_0^{2\pi}
    \innerprod{X(z),Y'(z)}\,dz. 
  \end{equation*}
  Then $\tilde L\mathfrak{g}=\reals\oplus L\mathfrak{g}$ with bracket
  \begin{equation*}
    [(a,X),(b,Y)] = (\omega(X,Y), [X,Y]).
  \end{equation*}
\item The extension $0\to \reals\to \tilde L\mathfrak{g}\to
  L\mathfrak{g}\to 0$ with bracket given by the inner product
  $\innerprod{,}$ on $\mathfrak{g}$ corresponds to an extension of Lie
  groups $1\to S^1\to \tilde LG\to LG\to 1$ if and only if
  $\innerprod{h_\alpha,h_\alpha}\in 2\integers$ for every coroot
  $h_\alpha$ of $G$.
\item If this integrality condition is satisfied, the extension
  $\tilde LG$ is uniquely determined. Moreover, there is a unique
  action of $\Diffeo^+(S^1)$ on $\tilde LG$ which covers the action on
  $LG$. In particular, there is a induced canonical lift of the action
  of $S^1$ on $LG$ by rotation of the argument to $\tilde LG$.
\item The integrality condition is satisfied if and only if
  $\omega/2\pi$, considered as an invariant differential form on 
  $LG$ (which is closed by the invariance of $\innerprod{,}$ and
  therefore of $\omega$), lifts to an integral cohomology class. It
  then represents the first Chern class of the principle $S^1$-bundle
  $S^1\to \tilde LG\to LG$ over $LG$. It follows that the topological
  structure completely determines the group extension.
\item If $G$ is simple and simply connected there is a
  \emph{universal} central extension $1\to S^1\to \tilde LG\to LG\to
  1$ (universal means that there is a unique map of extensions to any
  other central extension of $LG$).

  If $G$ is simple, all invariant bilinear forms on $\mathfrak{g}$ are
  proportional, and the universal extension corresponds to the
  smallest non-trivial one which satisfies the integrality
  condition. We call this also the basic inner product and the basic
  central extension.
\item For $SU(n)$ and the other simply laced groups, the basic inner
  product is the canonical inner product such that
  $\innerprod{h_\alpha,h_\alpha}=2$ for every coroot $\alpha$.
\item There is a precise formula for the adjoint and coadjoint action
  of suitable elements, compare Lemma \ref{lem:action of Waff on weights}.
\end{enumerate}

\section{Representations of loop groups}
\label{sec:repr-loop-groups}

\begin{definition}
  A representation of a loop group $LG$ (or more generally any
  topological group) is for us given by a locally convex topological
  vector space $V$ (over $\complexs$) with an action
  \begin{equation*}
    G\times V\to V; (g,v)\to gv
  \end{equation*}
  which is continuous and linear in the second variable.

  Two representations $V_1,V_2$ are called \emph{essentially
    equivalent}, if they contain dense $G$-invariant subspaces
  $V_1'\subset V_1$, $V_2'\subset V_2$ with a continuous
  $G$-equivariant bijection $V_1'\to V_2'$

  The representation $V$ is called \emph{smooth}, if there is a dense
  subspace of vectors $v\in V$ such that the map $G\to V;\,g\mapsto
  gv$ is smooth ---such vectors are called \emph{smooth vectors} of the
  representation.

  A representation $V$ is called irreducible if it has no closed
  invariant subspaces.
\end{definition}

\begin{example}
  The actions of $S^1$ on $C^\infty(S^1)$ ,$C^(S^1)$ and $L^2(S^1)$ by
  rotation of the argument are all equivalent and smooth.

  However, rotation of the argument does not define an action in our
  sense of $S^1$ on $L^\infty(S^1)$ because the corresponding map
  $S^1\times L^\infty(S^1)\to L^\infty(S^1)$ is not continuous.
\end{example}

\begin{remark}
  Given a representation $V$ of $LG$ and $z_0\in S^1$ (or more
  generally any diffeomorphism of $S^1$, $z_0$ gives rise to a
  diffeomorphism by translation), then we define a new representation
  $\phi^*V$ by composition with the induced automorphism of $LG$.

  We are most interested in representations which are symmetric,
  meaning that $\phi^*V\iso V$. Actually, we require a somewhat
  stronger condition in the following Definition \ref{def:proj_rep}.
\end{remark}

\begin{definition}\label{def:proj_rep}
  When we consider representation $V$ of $LG$, we really want to consider actions
  of $LG\semiProd S^1$, i.e.~we want an action of $S^1$ on $V$ which
  intertwines the action of $LG$; for $z_0\in S^1$ we want operators
  $R_{z_0}$ on $V$ such that $R_{z_0}fR_{z_0}^{-1}v = f(\cdot +z_0)v$
  for all $v\in V$, $f\in LG$.

  Moreover, we will study \emph{projective} representations,
  i.e.~representations such that $f_1\cdot (f_2\cdot v) = c(f_1,f_2)
  (f_1f_2)\cdot v$ with $c(f_1,f_2)\in \complexs\setminus\{0\}$.

  More precisely, these are actions of a central extension $1\to
  \complexs^{\times} \to \tilde LG\to LG\to 1$.

  Since $S^1$ acts on $V$, we get a decomposition
  $V=\overline{\bigoplus_{k\in\integers} V(k)}$, where $z\in S^1$ acts
  on $V(k)$ by multiplication with $z^{-k}$.

  We say $V$ is a representation of \emph{positive energy}, if
  $V(k)=0$ for $k<0$.
\end{definition}

\begin{example}
  The adjoint representation of $LG$ on $L\mathfrak{g}$, or the
  canonical representation of $LSU(n)$ on the Hilbert space
  $L^2(S^1,\complexs^n)$ are \emph{not} of positive energy, nor of
  negative energy. However, they are not irreducible.
\end{example}

\begin{remark}
  One can always modify the action of $S^1$ on a representation $V$ of
  $LG$ by multiplication with a character of $S^1$, so that ``positive
  energy'' and ``energy bounded below'' are more or less the same.

   The complex conjugate of a representation of negative energy is a
   representation of positive energy.
\end{remark}

\begin{proposition}
  An irreducible unitary representation $V$ of $\tilde LG\semiProd S^1$ which
  is of positive energy (and this makes only sense with the action of
  the extra circle) is also irreducible as a representation of $\tilde
  LG$.
\end{proposition}
\begin{proof}
  Let $T$ be the projection onto a proper $\tilde LG$-invariant summand of
  $V$. This operator commutes with the action of $\tilde LG$. Let
  $R_z$ be the operator through which $z\in S^1$ acts on $V$. We
  define the bounded operators $T_q:=\int_{S^1} z^q
  R_zTR_z^{-1}\;dz$. They all commute with $\tilde LG$, and $T_q$ maps
  $V(k)$ to $V(k+q)$.

  Let $m$ be the lowest energy of $V$. Then $T_q(V(m))=0$ for all
  $q<0$. Since $V$ is irreducible, $V$ is generated as a
  representation of $LG\semiProd S^1$ by $V(m)$. Since $V(m)$ is
  $S^1$-invariant, $T_q(V)=0$ for $q<0$. Since $T_{-q}=T_q^*$, we even
  have $T_q=0$ for all $q\ne 0$. Now, the $T_q$ are the Fourier
  coefficients of the loop $z\mapsto R_zTR_z^{-1}$. It follows that
  this loop is constant, i.e.~that $T$ commutes also with the action
  of $S^1$. But since $V$ was irreducible, this implies that $T=0$.

\end{proof}

The representation of positive energy of a loop group behave very much
like the representation of a compact Lie group. This is reflected in
the following theorem.

\begin{theorem}
  Let $G$ be a compact Lie group. Let $V$ be a smooth representation
  of positive energy of $\tilde LG\semiprod S^1$. Then, upto essential
  equivalence:
  \begin{enumerate}
  \item If $V$ is non-trivial then it does not factor through an
    honest representation of $LG$, i.e.~is truly projective.
  \item $V$ is a discrete direct sum of irreducible representations
    (of positive energy)
  \item $V$ is unitary
  \item The representation extends to a representation of $\tilde
    LG\semiProd \Diffeo^+(S^1)$, where $\Diffeo^+(S^1)$ denotes the
    orientation preserving diffeomorphisms and contains $S^1$ (acting
    by translation).
  \item $V$ extends to a holomorphic projective representation of $LG_\complexs$.
  \end{enumerate}
\end{theorem}

Granted this theorem, it is of particular importance to classify the
irreducible representations of positive energy

\begin{definition}
  Let $T_0\times T\times S^1\subset \tilde LG\semiProd S^1$ be a
  ``maximal torus'', with $T_0 = S^1$ the kernel of the central
  extension, $T$ a maximal torus of $G$ and $S^1$ the rotation group.
  We can then refine the energy decomposition of any representation $V$ to
  a decomposition
\begin{equation*}
  V=\overline{\bigoplus_{n,\alpha,h\in \integers\times \hat T\times
      \integers} V_{n,\alpha,h}}
  \end{equation*}
  according to the characters of $T_0\times T\times S^1$. The
  characters which occur are called the \emph{weights} of $V$. $n$ is
  called the \emph{energy} and $h$ the \emph{level}.
\end{definition}

\begin{lemma}\label{lem:action of Waff on weights}
  The action of $\xi\in\Hom(S^1,T)$ on a weight $(n,\alpha,h)$ is
  given by
  \begin{equation*}
    \xi(n,\alpha,h) = (n+\alpha(\xi)+ h\abs{\xi}^2/2, \alpha+h\xi^*,h)
  \end{equation*}
  where we identify $\xi$ with $\xi'(1)\in \mathfrak{t}$.

  Moreover, the norm is obtained from the inner product on
  $L\mathfrak{g}$ which corresponds to the central extension $\tilde
  LG$, and $\xi^*\in \hat T$ is the image of $\xi$ under the map
  $\mathfrak{t}\to\mathfrak{t}^*$ defined by this inner product.
\end{lemma}

\begin{remark}
    Note that $T_0$ commutes with every element of $\tilde LG\semiProd
  S^1$. Consequently, each representation decomposes into
  subrepresentations with fixed level, and an irreducible
  representation has only one level $h$.

  The level is a measure for the ``projectivity'' of the
  representation; it factors through an honest representation of $LG$
  if and only if the level is $0$.

  The weights of a representation are permuted by the normalizer of
  $T_0\times T\times S^1$, hence by the affine Weyl group
  $W_{aff}=\Hom(S^1,T)\semiProd W$ (where $W$ is the Weyl group of $G$).
\end{remark}

\begin{definition}
  Given a root $(n,\alpha)$ of $LG$, we define the \emph{coroot}
  $(-n\abs{h_\alpha}^2/2,h_\alpha)\in\reals\oplus \mathfrak{t}\subset
  \reals\oplus \mathfrak{t}\oplus\reals$, where we use again the inner
  product for the central extension $\tilde LG$.
\end{definition}


Our irreducible representation have a number of additional important
properties.
\begin{theorem}
  Let $V$ be a smooth irreducible representation of $LG$ of positive
  energy (i.e.~we really take a representation of $\tilde LG\semiProd
  S^1$). 
  \begin{enumerate}
  \item Then $V$ is of finite type, i.e.~for each energy $n$ the
    subspace $V(n)$ is finite dimensional. In particular, each weight
    space $E_{h,\lambda,n}$ is finite dimensional.
  \item $V$ has a unique lowest weight $(h,\lambda,n)$ with
    $E_{h,\lambda,n}\ne 0$. \emph{Lowest weight} means by definition,
    that for each positive root $(\alpha,m)$ the character
    $(h,\lambda-\alpha,n-m)$ does not occur as a weight in $V$.

    This lowest weight is \emph{antidominant}, i.e.~for each positive
    coroot $(-m\abs{h_\alpha}^2/2,h_\alpha)$ we have
    $\innerprod{(h,\lambda,n),(-m\abs{h_\alpha}^2/2),h_\alpha,-0}=
    \lambda(h_\alpha)-hm\abs{h_\alpha}^2/2 \le 0$.

    Since we have in particular to consider the positive roots
    $(\alpha,0)$ and $(-\alpha,1)$ (for each positive root $\alpha$ of
    $G$), this is equivalent to 
    \begin{equation}\label{eq:antidom_ineq}
      -h\abs{h_\alpha}^2/2\le \lambda(h_\alpha)\le 0
    \end{equation}
    for each positive root $\alpha$ of $G$.
  \item There is a bijection between isomorphism class of irreducible
    representations of $\tilde LG\semiProd S^1$ as above and antidominant weights.
  \end{enumerate}
\end{theorem}

\begin{corollary}
  If the level $h=0$, only $\lambda=0$ satisfies Inequality
  \eqref{eq:antidom_ineq}. In other words, among the representations
  considered here, only the trivial
  representation is an honest representation of $LG$, all others are
  projective.

  For a given level $h$, there are only finite possible antidominant weights
  (with $n=0$), because the $h_\alpha$ generate $\mathfrak{t}$. 
\end{corollary}

\begin{example}
  If $G$ is simple and we look at an antidominant weight
  $(h,\lambda,0)$, then $-\lambda$ is a dominant weight in the usual
  sense of $G$, i.e.~contained in the corresponding simplicial cone in
  $\mathfrak{t}^*$, but with the extra condition that it is contained
  in the simplex cut off by $\{\mu\mid \mu(\alpha_0)=h\}$, where
  $\alpha_0$ is the highest weight of $G$.

  We get in particular the so called \emph{fundamental weights}
  \begin{enumerate}
  \item $(1,0,0)$
  \item $(\innerprod{\omega_i,\alpha_0}, -\omega_i,0)$, with
    $\omega_i$ the fundamental weights of $G$ determined by $\omega_i(h_{\alpha_j})=\delta_{ij}$.
  \end{enumerate}
  The antidominant weights are exactly the linear combinations of the
  fundamental weights with coefficients in $\naturals\cup\{0\}$.
\end{example}

Given an irreducible representation $V$ of $\tilde LG\semiProd S^1$ of
lowest weight $(h,\lambda,0)$,
to get a better understanding of the we want to determine which other
weights occur in $V$ (or rather, we want to find restrictions for
those weights).

First observation: the whole orbit under $W_{aff}$ occurs. This
produces, for the $\eta\in \Hom(S^1,T)$, the weights
$(h,\lambda+h\eta^*,\lambda(\eta)+h\abs{\eta}^2/2)$,

\begin{example}
  If $G=SU(2)$, we have the isomorphism $(\hat T\subset \mathfrak{t}
  )\iso (\integers\subset \reals)$, where $\left(
    \begin{smallmatrix}
      2\pi it & 0\\ 0 & -2\pi it
    \end{smallmatrix}\right)\in\mathfrak{t}$ is mapped to
  $t\in\reals$, $\mathfrak{t}^*$ is identified with $\mathfrak{t}$
  using the standard inner product and this way the character
  $\diag(z,z^{-1})\mapsto z^\mu$ in $\hat T$ is mapped to
    $\mu\in\integers$.

    Under this identification, the lowest weight $\alpha_0$ is
    identified with $1\in\integers$.

    There are exactly two fundamental weights. The $W_{aff}$-orbit of the weight
    $(1,0,0)$ is (with this identification of $\hat T$ with
    $\integers$)  $\{(1,2k,m)\mid (2k)^2= 2m\}$, the set of all
    weights of the corresponding irreducible representations is
    $\{(1,2k,m)\mid (2k)^2\le 2m\}$. Similarly, the orbit of
    $(1,-1,0)$ is  $\{(1,2k+1,m)\mid (2k+1)^2= 2m+1\}$, the set of
    all weights is $\{(1,2k+1,m)\mid (2k+1)^2\le 2m+1\}$.

    In general (for arbitrary $G$), the orbit of the lowest weight $(h,\mu,m)$ (with
    $\abs{(h,\mu,m)}^2=\abs{\mu}^2+2mh$) is contained in the parabola
    $\{(h,\mu',m')\mid \abs{\mu'}^2=\abs{(h,\mu,m)}^2+2m'h\}$.

    All other weights are contained in the interior of this parabola
    (i.e.~those points with ``$=$'' replaced by ``$\le$'').

    The last statement follows because, by translation with an element
    of $W_{aff}$ we can assume that $(h,\mu',m')$ is antidominant. Then
    \begin{equation*}
      \abs{(h,\mu',m')}^2- \abs{(h,\mu,m)}^2 =
      \innerprod{(h,\mu',m')+(h,\mu,m), (h,\mu',m')-(h,\mu,m)}\le 0,
    \end{equation*}
    because the first entry is antidominant and the second one is
    positive, $(h,\mu,m)$ being a lowest weight.

    We use the fact (not proved here)  that the extension inner product on $\mathfrak{t}$ extends
    to an inner product on $\reals\oplus\mathfrak{t}\oplus\reals$
    which implements the pairing between roots and coroots.
\end{example}

\section{Proof for Section \ref{sec:centr-extens-lg} and Homogeneous spaces of $LG$}
\label{sec:homog-spac-lg}

We now want to indicate the proofs of the statements of Section
\ref{sec:repr-loop-groups}. 

The basic idea is that we can mimic the we can mimic the Borel-Weil
theorem for compact Lie groups. It can be stated as follows:

\begin{theorem}
  The homogeneous space $G/T$ has a complex structure, because it is
  isomorphic to $G_\complexs/B^+$, where $G_\complexs$ is the
  complexification of $G$ and $B^+$ is the Borel subgroup. In case
  $G=U(n)$, $G_\complexs=Gl(n,\complexs)$ and $B^+\subset
  Gl(n,\complexs)$ is the subgroup of upper triangular matrices; the
  homogeneous space is the flag variety.

  To each weight $\lambda\colon T\to S^1$ there is a uniquely
  associated holomorphic line bundle $L_\lambda$ over
  $G_\complexs/B^+$ with action of $G_\complexs$.

  $L_\lambda$ has non-trivial holomorphic sections if and only if
  $\lambda$ is an antidominant weight. In this case, the space of
  holomorphic section is an irreducible representation with lowest
  weight $\lambda$.
\end{theorem}

For loop groups, the relevant homogeneous space is
\begin{equation*}
  Y:= LG/T = LG_{\complexs}/B^+G_\complexs,\qquad B^+G_\complexs =\{
  \sum_{k\ge 0} \lambda_kz^k\mid \lambda_0\in B^+\}
  \end{equation*}
  Recall that for $GL(n,\complexs)$, the \emph{Borel subgroup} $B^+$ is the
  subgroup of upper triangular matrices.

Note that the second description defines on $Y$ the structure of a
complex manifold.

  We have also $Y=\tilde LG/\tilde T = \tilde LG_\complexs/ \tilde B^+G_\complexs$.

  \begin{lemma}
    Each character $\lambda\colon \tilde T\to S^1$ has a unique
    extension $\tilde B^+G_\complexs = \tilde T_\complexs\cdot \tilde
    N^+G_\complexs \to \complexs^x$,, with
    $N^+G_\complexs:=\{\sum_{k\ge 0} \lambda_k z^k\mid \lambda_0\in
    N^+G_\complexs\}$, where $N^+$ is the nilpotent subgroup of
    $G_\complexs$ whose Lie algebra is generated by the positive root
    vectors, for $Gl(n,\complexs)$ it is the group of upper triangular
    matrices with $1$s on the diagonal. This way defines a
    holomorphic line bundle 
    \begin{equation*}
L_\lambda:= \tilde LG_\complexs\times_{\tilde
  B^+G_\complexs}\complexs\text{ over }Y.
\end{equation*}

    Write $\Gamma_\lambda$ for the space of holomorphic sections of
    $L\lambda$. This is a representation of $\tilde LG\semiProd S^1$.
  \end{lemma}

  \begin{lemma}
    The space $Y=LG/T$ contains the affine Weyl group
    $W_{aff}=(\Hom(S^1,T)\cdot N(T))/T$, where $N(T)$ is the
    normalizer of $T$ in $G\subset LG$.

    $Y$ is stratified by the orbits of $W_{aff}$ under the action of
    $N^-LG_\complexs:= \{\sum_{k\le 0} \lambda_kz^k\mid \lambda_0\in
    N^-G_\complexs\}$. Here $N^-G_\complexs$ is the nilpotent Lie
    subgroup whose Lie algebra is spanned by the negative root vectors
    of $\mathfrak{g}_\complexs$. For $Gl(n,\complexs)$ this is the
    group of lower triangular matrices with $1$s on the diagonal.
  \end{lemma}

  \begin{theorem}
    Assume that $\lambda\in \Hom(\tilde T,S^1)$ is a weight such that
    the space $\Gamma_\lambda$ of holomorphic sections of $L_\lambda$
    is non-trivial. Then
    \begin{enumerate}
    \item $\Gamma_\lambda$ is a representation of positive energy.
    \item $\Gamma_\lambda$ is of finite type, i.e.~each fixed energy subspace
      $\Gamma_\lambda(n)$ is finite dimensional
    \item $\lambda$ is the lowest weight of $\Gamma_\lambda$ and is antidominant.
    \item $\Gamma_\lambda$ is irreducible.
    \end{enumerate}
  \end{theorem}
  \begin{proof}
    We use the stratification of $Y$ to reduce to the top stratum. It
    turns out that the top stratum is $N^-LG_\complexs$ and that
    $L_\gamma$ trivializes here, so $\Gamma_\lambda$ restricts to the
    space of holomorphic functions on $N^-LG_\complexs$. Because
    holomorphic sections are determined by their values on the top
    stratum, this restriction map is injective. We can
    further, with the exponential map (which is surjective in this
    case), pull back to holomorphic functions on
    $N^-L\mathfrak{g}_\complexs$, and then look at the Taylor
    coefficients at $0$.

    This way we finally map invectively into $\prod_{p\ge 0}
    S^p(N^-L\mathfrak{g}_\complexs)^*$, where $S^p(V)^*$ is the space
    of $p$-multilinear maps $V\times\cdots\times V\to \complexs$. This
    map is indeed $\tilde T\times S^1$-equivariant, if we multiply the
    obvious action on the target with $\lambda$.

    It now turns out that $N^-L\mathfrak{g}_\complexs$ has (essentially
    by definition) negative energy, and therefore the duals
    $S^p(N^-L\mathfrak{g}_\complexs)$ all have positive
    energy; and the weights are exactly the positive roots (before
    multiplication with $\lambda$). Consequently $\lambda$ is of
    lowest weight. If for some positive root $(\alpha,n)$, we had
    $\lambda(\alpha,n)=m>0$, then reflection in $W_{aff}$
    corresponding to $\alpha$ would map $\lambda$ to $\lambda
    -m\alpha$, which on the other hand can not be a root of
    $\Gamma_\lambda$ if $\lambda$ is a root. Consequently, $\lambda$ is
    antidominant.

   Explicit calculations also show that the image of the
    ``restriction map'' is contained in a subspace of finite type.

    To prove that $\Gamma_\lambda$ is irreducible, we look at the
    subspace of lowest energy. This is a representation of
    $G_\complexs$. Pick a lowest weight vector for this
    representation, it is then invariant under the nilpotent subgroup
    $N^-$. Since it is of lowest energy, it is even invariant under
    $N^-LG_\complexs$. On the other hand, since the top stratum of $Y$
    is $N^-LG_\complexs$, the value of an invariant section at one
    point completely determines it, so that the space of such sections
    is $1$-dimensional. $B^-LG_\complexs$ acts on this space by
    multiplication with the holomorphic homomorphism $\lambda\colon
    B^-LG_\complexs\to\complexs^\times$, so that $\lambda$ really
    occurs as lowest weight.

    We now show that this vector is actually a cyclic vector,
    i.e.~generates $\Gamma_\lambda$ under the action of $\tilde
    LG\semiProd S^1$. Else choose a vector of lowest energy not in
    this subrepresentation, of lowest weight for the corresponding
    action of the compact group $G$, and we get with the argument as
    above a second $N^-LG_\complexs$ invariant section.
  \end{proof}

  \begin{lemma}
 $\Gamma_\lambda\ne 0$ if and only   If $\lambda$ is antidominant.
  \end{lemma}
  \begin{proof}
    If $\alpha$ is a positive root of $G$ ($(\alpha,0$ therefore a
    positive root of $LG\semiProd S^1$), we get a corresponding
    inclusion $i_\alpha\colon Sl_2(\complexs)\to \tilde LG_\complexs$
    whose restriction to $\complexs^\times\subset Sl_2(\complexs)$ is
    the exponential of $h_\alpha$. Since $\alpha$ is positive,
    $B^+Sl_2(\complexs)$ is mapped to $B^+\tilde
    LG_\complexs$. Therefore, we get an induced map
    $P^1(\complexs)=Sl_2(\complexs)/B^+ \to Y$. The pullback of $L_\lambda$ under
    this map is the line bundle associated to $\lambda\circ
    h_\alpha$. If $\Gamma_\lambda$ is non-trivial, we can therefore
    pull back to obtain a non-trivial holomorphic section of this
    bundle over $P^1(\complexs)$. These exist only if $\lambda\circ
    h_\alpha=\lambda(h_\alpha)\le 0$.

    This prove half the conditions for antidominance. We omit the
    prove of the other half, where we have to consider the positive
    root $(-\alpha,1)$.

    If, on the other hand, $\lambda$ is antidominant, one constructs a
    holomorphic section along the stratification of $Y$.
  \end{proof}

  \begin{theorem}
    An arbitrary smooth representation of $\tilde LG\semiProd S^1$ of
    positive energy splits (upto essential equivalence) as a direct
    sum of representations of the form $\Gamma_\lambda$.

    In particular, the $\Gamma_\lambda$ are exactly the irreducible representations.
  \end{theorem}
  \begin{proof}
    If $E$ is a representation of positive energy, so is $\overline
    E^*$. Pick in the $G$-representation $\overline E^*(0)$ a vector
    $\epsilon$ of lowest weight $\lambda$. For each smooth vector
    $v\in E$, the map 
    \begin{equation*}
      s_v\colon \tilde LG_\complexs\to \complexs; f\mapsto \epsilon(f^{-1}v)
    \end{equation*}
    turns out to define a holomorphic section of $L_\lambda$. This
    give a non-trivial map $E\to \Gamma_\lambda$. If $E$ is irreducible it
    therefore is essentially equivalent to $\Gamma_\lambda$ 

    Using $\overline\Gamma_\lambda^*$ and similar constructions, we
    can split off factor $\Gamma_\lambda$ successively from an
    arbitrary representation $E$.
  \end{proof}

The proofs of the statements about the structure of these homogenous
spaces uses the 
study of related Grassmannians. These we define and study for $LU(n)$;
results for arbitrary compact Lie groups follow by embedding into
$U(n)$ and reduction to the established case.

Let $H=H_+\oplus H_-$ be a (polarized) Hilbert space. The important
example for us is $H=L^2(S^1,\complexs^n)$, with $H_+$ generated by
functions $z^k$ for $k\ge 0$ and $H_-$ generated by $z^k$ with $k<0$
(negative or positive Fourier coefficients vanish).

On this Hilbert space, the complex loop group $LGl(n,\complexs)$ acts
by pointwise multiplication.

We define the Grassmannian $Gr(H):=\{W\subset H\mid \pr_+\colon W\to
H_+\text{ is Fredholm}, pr_=\colon W\to H_-\text{ is
  Hilbert-Schmidt}\}$. This is a Hilbert manifold.

We define the restricted linear group $Gl_{res}(H):= \{\left(
  \begin{smallmatrix}
    a & b\\c & d
  \end{smallmatrix}\right)\mid b,c\text{ Fredholm}\}$. (This implies
that $a,d$ are Hilbert-Schmidt). Set $U_{res}(H):=U(H)\cap
Gl_{res}(H)$. These groups act on $Gr(H)$. 

\begin{lemma}
  The image of $LGl(n,\complexs)$ in $Gl(H)$ is contained in $Gl_{res}(H)$.
\end{lemma}
\begin{proof}
  Write $\gamma=\sum_k{\gamma_k} e^{ik\theta} =
  \begin{pmatrix}
    a & b\\ c & d
  \end{pmatrix}\in LGl(n,\complexs)$. Then for each $n>0$,
  $b(e^{-in\theta}) =\sum_{k\ge n }\gamma_k e^{i(k-n)\theta}$, and for
  $n\ge 0$, $c(e^{in\theta})=\sum_{k<-n} \gamma_k
  e^{i(k+n)\theta}$. Therefore $\norm{b}_{HS} = \sum_{n\ge
    0}\sum_{k\ge n} \abs{\gamma_k}^2 = \sum_{k\ge 0}
  (1+k)\abs{\gamma_k}^2<\infty$ since smooth functions have rapidly
  decreasing Fourier expansion. Similarly for $c$.
\end{proof}

$Gr(H)$ has a cell decomposition/stratification in analogy to the
finite dimensional Grassmannians (with infinitely many cells of all
kinds of relative dimensions). Details are omitted here, but this is
an important tool in many proofs.

$Gr(H)$ contains the subspace $Gr_\infty(H):=\{W\in Gr(H)\mid
\im(pr_-)\cup \im(\pr_+)\subset C^\infty(S^1,\complexs^n)\}$. Inside
this one we consider the subspace $Gr_\infty^{(n)}:=\{W\in
Gr_\infty(H)\mid zW\subset W\}$.

It turns out that $Gr_\infty(H)^{(n)}= LGl(n,complexs)/
L^+Gl(n,\complexs)= LU(n)/U(n)$.

The main point of the definition of $Gr(H)$ is that we can define a
useful (and fine) virtual dimension for its elements; measuring the dimension of
$W\cap z^mH_-$ for every $m$. This is used to stratify these
Grassmannians and related homogeneous spaces, and these
stratifications were used in the study of the holomorphic sections of
line bundles over these spaces.

\nocite{MR900587,MR1853241}
{\small
\bibliographystyle{plain}
\bibliography{Loop_groups}
}

\end{document}